\documentstyle[12pt]{article}
\def\llim{\mathop{\rm \underline {lim}}}

\def\b{{\beta}}

\def\e{{\varepsilon}}

\def\R{{{\bf R}^2}}
\def\D{{{\cal D}}}
\def\g{{\gamma}}
\def\v{{\vert}}
\def\A{{{\cal A}_\rho(\D)}}
\def\ee{{\rm exp}}
\def\gr{{\rm grad}}
\def\Im{{\rm Im}\,}
\def\ch{{\rm ch}\,}
\def\sh{{\rm sh}\,}
\def\F{{\Phi(y,x)}}

\date{14 February 1995}
\author{Z.~R.~Ashurova, \  Y.~I.~Zhuraev}
\title{Uniqueness theorem for unbounded domain}

\begin{document}
\maketitle
\begin{abstract}
We prove in this paper the uniqueness theorem for a certain class of
harmonic functions defined in unbounded domain lying in a band.

\end{abstract}

\section{Introduction}
Let $\D$ be a simply connected unbounded domain in $\R$ lying inside a band
of minimal width defined by conditions $0\leq y_2\leq h$,
$h=\frac{\pi}{\rho}$, $\rho >0$, $y=(y_1,y_2)$. We assume that $\partial\D=
\g_1\cup\g_2$, where $\g_i$ are smooth curves given by formulas
$y_2=f_i(y_1),\ i=1,2$ with $$ \v f_i(y_1)\v\leq {\rm Const},\quad
\v f'_i(y_1)\v\leq {\rm Const}. $$
Let $\A$ be a space of harmonic functions in $\D$ which are continuous with
their partial derivatives of the first order up to the ending points of
$\partial\D$ and satisfying the following condition:
$$ \v U(y)\v +\v\gr U(y)\v\leq \ee({\rm o}(\ee\rho\v y_1\v)),\ y_1\to\infty.
$$
The purpose of this paper is to prove the following

{\bf Theorem.}\ {\it If function $U\in\A$ satisfies the two following
conditions
\begin{equation}\label{1}
U(y)=0, \quad y\in\partial \D,
\end{equation}
\begin{equation}\label{2}
\left|\frac{\partial U}{\partial{\bf n}}(y)\right|\leq\ee(c\v y\v),\
y\in\partial\D,\ c={\rm Const}<\frac{\rho}{2},
\end{equation}
(where ${\bf n}$ is the exterior normal vector for the boundary
$\partial\D$) then $U\equiv 0$ in $\D$.}

In the case when the condition (\ref{2}) is changed by
$\left|\frac{\partial U}{\partial{\bf n}}(y)\right|\leq c_1\v y\v^k$ with
$c_1={\rm Const}$, $y\in\partial\D$, and $k$ being an non-negative integer
this theorem has been proved by Z.~R.~Ashurova in~\cite{1}.
The proof of this theorem is based on integral representation of harmonic
functions in unbounded domain with non-compact boundary~\cite{3} and on the
theorem 6.1 by E.~M.~Landis~\cite{2}.

\section{Crucial inequality}
Let $y=(y_1,y_2)$, $x=(x_1,x_2)$, $x,y\in\R$, $0<x_2,y_2<h$,
$h=\frac{\pi}{\rho}$ with $\rho>0$. For $x\not= y$ we can define a function
$\F$ by formula
\begin{equation}\label{3}
\F=-\frac{1}{2\pi K(x_2)}\int_0^\infty \Im\frac{K(y_2+i\eta)}{y_2-x_2+i\eta}
\cdot\frac{u\,du}{\eta},
\end{equation}
where $\eta^2=u^2+\alpha^2$, $\alpha^2=(y_1-x_1)^2$,
\begin{equation}\label{4}
K(\omega)=(\omega+3h-x_2)^{-1}\ee\left(-a\,\ch
i\rho_1(\omega-\frac{h}{2})\right), \qquad \omega=\xi+i\eta,
\end{equation}
 $a$, $\rho_1$ are positive numbers and $0<\rho_1<\rho$.
Then for the function $\F$ the estimate
\begin{equation}\label{5}
\left|\F\right|\leq\frac{C_0}{\ee\left(a\,\cos(\rho_1(y_2-\frac{h}{2}))\,\ch
\rho_1\alpha\right)}\left(1+\ln\left(1+\frac{15h^2}{\alpha^2+\beta^2}\right)
\right),
\end{equation}
is valid (here $\beta^2=(y_2-x_2)^2$ and $C_0={\rm Const}$). The aim of
this section is to prove the inequality (\ref{5}).

If $\omega=x_2$ then by (\ref{4}) we have $$
K(x_2)=(3h)^{-1}\ee\left(a\,\cos \rho_1(x_2-\frac{h}{2})\right). $$
As we have $-\frac{\pi}{2}<\rho_1(x_2-\frac{h}{2})<\frac{\pi}{2}$, so
$\cos\rho_1(x_2-\frac{h}{2})>0$, henceforth
\begin{equation}\label{6}
\frac{1}{K(x_2)}=3h\,\ee\left(a\,\cos\rho_1(x_2-\frac{h}{2})\right)\leq
3h\,\ee(a).
\end{equation}
Put $\beta_1=y_2-x_2+3h$, $r_1^2=\alpha^2+\beta_1^2$,
$r^2=\alpha^2+\beta^2$. It can be easily seen that
\begin{eqnarray}\label{7}
\Im\frac{K(y_2+i\eta)}{(y_2-x_2+i\eta)}&=&
\frac{(\b\b_1-\eta^2)\sin a\left
(\sin
\rho_1(y_2-\frac{h}{2})\,\sh \rho_1\eta\right)}
{(u^2+r^2)(u^2+r_1^2)\ee\left(a\cos \rho_1(y_2-\frac{h}{2})\,\ch \rho_1\eta
\right)}
\nonumber \\
& -&
\frac{\eta(\b+\b_1)\cos \left(a\sin
\rho_1(y_2-\frac{h}{2})\,\sh \rho_1\eta\right)}
{(u^2+r^2)(u^2+r_1^2)\ee\left(a\cos \rho_1(y_2-\frac{h}{2})\,\ch \rho_1\eta
\right)}.
\end{eqnarray}
As we have $\alpha\leq\eta$, so
$$ \ee\left(a\cos \rho_1(y_2-\frac{h}{2})\,\ch \rho_1\alpha\right)\leq
 \ee\left(a\cos \rho_1(y_2-\frac{h}{2})\,\ch \rho_1\eta\right). $$
Since
$$ \lim_{\eta\to 0}\frac{\b\b_1\sin(a\,\sh\rho_1\eta)}{\eta^2+\b^2}=0, $$
there exists such $\e>0$ that for $\v\eta\v<\e$ the function
$$\frac{\b\b_1\sin(a\,\sh\rho_1\eta)}{\eta^2+\b^2} $$ is bounded. If
$\v\eta\v\geq\e$ then we have the inequality
$$ \left| \frac{\b\b_1\sin(a\,\sh\rho_1\eta)}{\eta^2+\b^2}\right|\leq\frac
{4h}{\e^2}. $$
Moreover, we have $\frac{\eta^2}{u^2+r^2}=\frac{\eta^2}{\eta^2+\b^2}\leq 1$.
Thus we obtain
\begin{eqnarray}\label{8}
\left|
\frac{(\b\b_1-\eta^2)\sin a\left
(\sin
\rho_1(y_2-\frac{h}{2})\,\sh \rho_1\eta\right)}
{(u^2+r^2)\ee\left(a\cos \rho_1(y_2-\frac{h}{2})\,\ch \rho_1\eta
\right)} \right|
\nonumber \\
\leq \frac{C_1}{\ee\left(a\cos \rho_1(y_2-\frac{h}{2})\,\ch \rho_1\alpha
\right)},
\end{eqnarray}
where $C_1$ is some constant. (Further by $C_i$ we shall denote conctants)
It is not difficult to obtain the following estimate:
\begin{eqnarray}\label{9}
\left|\frac{\eta(\b+\b_1)\cos \left(a\sin
\rho_1(y_2-\frac{h}{2})\,\sh \rho_1\eta\right)}
{\ee\left(a\cos \rho_1(y_2-\frac{h}{2})\,\ch \rho_1\eta
\right)}\right| \nonumber \\
\leq \frac{C_2h}{\ee\left(a\cos \rho_1(y_2-\frac{h}{2})\,\ch \rho_1\alpha
\right)}.
\end{eqnarray}
Estimate fromabove the function $\F$. By (\ref{6}) -- (\ref{9}) we obtain
\begin{eqnarray*}
\v\F\v&\leq&
\frac{C_3}{\ee\left(a\cos \rho_1(y_2-\frac{h}{2})\,\ch \rho_1\alpha\right)}
\\
&\times&\left(\int_0^\infty\frac{u\,du}{(u^2+r_1^2)\sqrt{u^2+r^2}}+
\int_0^\infty\frac{u\,du}{(u^2+r^2)(u^2+r_1^2}\right).
\end{eqnarray*}
As $$\frac{u}{\sqrt{u^2+\alpha^2}}\leq 1,$$ the integral
$$\int_0^\infty\frac{du}{u^2+r_1^2}$$ is convergent and as
$$\int_0^\infty\frac{u\,du}{(u^2+r^2)(u^2+r_1^2}=\ln(1+\frac{15h^2}{r^2}),$$
so $$
\left|\F\right|\leq\frac{C_0}{\ee\left(a\,\cos(\rho_1(y_2-\frac{h}{2}))\,\ch
\rho_1\alpha\right)}\left(1+\ln\left(1+\frac{15h^2}{\alpha^2+\beta^2}\right)
\right). $$
Thus (\ref{5}) is proved.

\section{Proof of theorem}
Provided that the conditions of the theorem are valid we have~\cite{3}
$$
\int_{\partial\D}\left(U\frac{\partial\Phi}{\partial{\bf n}}-\Phi
\frac{\partial U}{\partial{\bf n}}\right)ds=
\left\lbrace
\begin{array}{cc}
0,& x\notin \D\cup\partial\D \\
U(x), & x\in \D.
\end{array}
\right.
$$
Using (\ref{1}) this formula can be written in the form
\begin{equation}\label{10}
U(x)=-\int_{\g_1}\F\frac{\partial U}{\partial{\bf n}}ds-
\int_{\g_2}\F\frac{\partial U}{\partial{\bf n}}ds, \quad x\in\D.
\end{equation}
Denote
$$ I_j=-\int_{\g_j}\F\frac{\partial U}{\partial{\bf n}}ds,\quad j=1,2,\quad
a_1=a\cos\rho_1\frac{h}{2}. $$
Since $\g_1$ is given by the equation $y_2=f_1(y_1)$ where $f_1$ is a
bounded function with abounded derivative, by (\ref{2}) and (\ref{5}) we
obtain
\begin{eqnarray*}
\v I_1\v&\leq& C_4\int_{-\infty}^\infty\frac{\ee(c\sqrt{y_1^2+h^2})}{\ee(a_1\,
\ch \rho_1\alpha)}\left(1+\ln(1+\frac{15h^2}{r^2})\right)\,dy_1 \\
&=&C_4(I_{11}+I_{12}),
\end{eqnarray*}
where we denote by $I_{11}$ (resp. by $I_{12}$) the integral on the
interval $(-\infty,0)$ (resp. $(0,\infty)$). Estimate at first $I_{12}$.
We can write it in the form
\begin{eqnarray}\label{11}
I_{12}
&=&\int_0^\infty\frac{\ee(c\sqrt{y_1^2+h^2})}{\ee(a_1\,\ch\rho_1\alpha)}dy
_1 \nonumber \\
&+&\int_0^\infty\ln\left(1+\frac{15h^2}{r^2}\right)
\frac{\ee(c\sqrt{y_1^2+h^2})}{\ee(a_1\,\ch\rho_1\alpha)}dy_1
\end{eqnarray}
As $\sqrt{y_1^2+h^2}\leq\sqrt{y_1^2+h^2+2y_1h}=y_1+h$
when$y_1\in(0,+\infty)$, so $$
\int_0^\infty\frac{\ee(c\sqrt{y_1^2+h^2})}{\ee(a_1\,\ch\rho_1\alpha)}dy_1
\leq
\int_0^\infty\frac{\ee c(y_1+h)}{\ee(a_1\,\ch\rho_1(y_1-x_1))}dy_1. $$
By change of variables $t=y_1-x_1$
we can write the last integral in the form
\begin{eqnarray}\label{12}
\int_0^\infty\frac{\ee c(y_1+h)}{\ee(a_1\,\ch\rho_1(y_1-x_1))}dy_1
&=&
\ee c(x_1+h)
\int_{-x_1}^\infty\frac{\ee (ct)}{\ee(a_1\,\ch\rho_1t)}dy_1 \nonumber \\
\leq C_5\ee(cx_1)\cdot
\int_{-\infty}^\infty\frac{\ee (ct)}{\ee(a_1\,\ch\rho_1t)}dy_1 &\leq &
C_6\ee(cx_1),
\end{eqnarray}
because the integral $$
\int_{-\infty}^\infty\frac{\ee (ct)}{\ee(a_1\,\ch\rho_1t)}dy_1 $$ is
convergent.
One can show that for the second integral a similar estimate is valid
\begin{equation}\label{13}
\left|\int_0^\infty\ln\left(1+\frac{15h^2}{r^2}\right)
\frac{\ee(c\sqrt{y_1^2+h^2})}{\ee(a_1\,\ch\rho_1\alpha)}dy_1\right|\leq
C_7\ee(cx_1).
\end{equation}
In the way similar to $I_{12}$ we have
\begin{equation}\label{14}
\v I_{11}\v\leq C_8\ee(cx_1).
\end{equation}
By inequalities (\ref{12}) -- (\ref{14}) we have
\begin{equation}\label{15}
\left|\int_{\g_1}\F\frac{\partial U}{\partial {\bf n}}ds\right|
\leq C_9\ee(cx_1)\leq C_9\ee(c\v x\v).
\end{equation}
The estimate of our integral on the path $\g_2$ can be made in a similar
way.
Thus by (\ref{10}) and (\ref{15}) we have $$ \v U(x)\v\leq C_{10}\ee(c\v
x\v). $$
As $c<\frac{\rho}{2}=\frac{\pi}{2h}$, so $$ \llim_{R\to\infty}\frac
{\max_{\v x\v=R}\v U(x)\v}
{\ee\frac{\pi}{2h}\v x \v}=0. $$
Hence using the theorem 6.1 of~\cite{2} we obtain that $U(x)\leq 0$.
Applying the above arguments for the function $-U(x)$ we obtain that
$-U(x)\leq 0$. Thus $U(x)\equiv 0$, $x\in\D$. We are done.

{\bf Acknowledgement.}\
The second author thanks for partial support by
the International Science Foundation (Soros) (grant N MGM000). We are
grateful to Sh.~Ya.~Yarmuhamedov and A.~S.~Mishchenko for helpful discussions.

\vspace{2.5cm}
\noindent
Ashurova Z.~R.,\quad Zhuraev Y.~I.\\
Dept. of Mathematics \\
Samarkand State University \\
Samarkand, 703004, Uzbekistan
\end{document}